\documentclass{amsart}
\usepackage{amsmath,amssymb}
\newtheorem{theorem}{Theorem}[section]
\newtheorem{lemma}[theorem]{Lemma}
\newtheorem{proposition}[theorem]{Proposition}
\newtheorem{corollary}[theorem]{Corollary}

\newtheorem{remark}[theorem]{Remark}

\begin{document}
\title{Flat modules over valuation rings}
\author{Fran\c{c}ois Couchot}
\address{Laboratoire de Math\'ematiques Nicolas Oresme, CNRS UMR
  6139,
D\'epartement de math\'ematiques et m\'ecanique,
14032 Caen cedex, France}
\email{couchot@math.unicaen.fr} 



\begin{abstract} Let $R$ be a valuation ring and let $Q$ be its total quotient ring. It is proved that any singly projective (respectively flat) module is finitely projective if and only if $Q$ is maximal (respectively artinian). It is shown that each singly projective module is a content module if and only if any non-unit of $R$ is a zero-divisor and that each singly projective module is locally projective if and only if $R$ is self injective. Moreover, $R$ is maximal if and only if each singly projective module is separable, if and only if any flat content module is locally projective. Necessary and sufficient conditions are given for a valuation ring with non-zero zero-divisors to be strongly coherent or $\pi$-coherent.

A complete characterization of semihereditary commutative rings which are $\pi$-coherent is given. When $R$ is a commutative ring with a self FP-injective quotient ring $Q$, it is  proved that  each flat $R$-module is finitely projective if and only if $Q$ is perfect.
\end{abstract}
\maketitle

In this paper, we consider the following properties of modules: P-flatness, flatness, content flatness, local projectivity, finite projectivity and single projectivity. We investigate the relations between these properties when $R$ is a PP-ring or a valuation ring. Garfinkel (\cite{Gar76}),  Zimmermann-Huisgen (\cite{ZH76}), and Gruson and Raynaud (\cite{GrRa71}) introduced the concepts of locally projective modules and strongly coherent rings and developed important theories on these. The notions of finitely projective modules and $\pi$-coherent rings are due to Jones (\cite{Jon82}). An interesting study of finitely projective modules and singly projective modules is also done by Azumaya in \cite{Azu87}.
For a module $M$ over a ring $R$, the following implications always hold:
\[\begin{matrix}
M\ \mathrm{is\ projective} & \Rightarrow & M\ \mathrm{is\ locally\ projective} & \Rightarrow & M\ \mathrm{is\ flat\ content} \\
{} & {} & \Downarrow & {} & \Downarrow \\
{} & {} & M\ \mathrm{is\ finitely\ projective} & \Rightarrow & M\ \mathrm{is\ flat} \\
{} & {} & \Downarrow & {} & \Downarrow \\
{}  & {} & M\ \mathrm{is\ singly\ projective} & \Rightarrow & M\ \mathrm{is\ P-flat},
\end{matrix}\]
but there are not generally reversible. However, if $R$ satisfies an additional condition, we get some equivalences. For instance, in \cite{Bas60}, Bass defined a ring $R$ to be right perfect if each flat right module is projective. In \cite{ZH80} it is proved that a ring $R$ is right perfect if and only if each flat right module is locally projective, and if and only if each locally projective right module is projective. If $R$ is a commutative arithmetic ring, i.e. a ring whose lattice of ideals is distributive, then any P-flat module is flat. By \cite[Proposition 16]{Azu87}, if $R$ is a commutative domain, each P-flat module is singly projective, and by \cite[Proposition 18 and 15]{Azu87} any flat left  module is finitely projective if $R$ is a commutative arithmetic domain or a left noetherian ring. Consequently, if $R$ is a valuation domain each P-flat module is finitely projective. When $R$ is a valuation ring, we prove that this result holds if and only if the ring $Q$ of quotients of $R$ is artinian. Moreover, we show that $R$ is maximal if and only if any singly projective module is separable or any flat content module is locally projective, and that $Q$ is maximal if and only if each singly projective module is finitely projective.

In Section~\ref{S:baer}, necessary and sufficient conditions are given for a commutative semihereditary ring to be $\pi$-coherent. Moreover we characterize commutative PP-rings for which each product of singly projective modules is singly projective.

In the last section  we study the valuation rings $R$ for which each product of content (respectively singly, finitely, locally projective) modules is content (respectively singly, finitely, locally projective). The results are similar to those obtained by Zimmermann-Huisgen and Franzen in \cite{Fra84}, and by Kemper in \cite{Kem00}, when $R$ is a domain. However, each valuation domain is $\pi$-coherent but not necessarily strongly coherent. We prove that a valuation ring with non-zero zero-divisors is $\pi$-coherent if and only if it is strongly coherent.

\section{Definitions and preliminaries}
\label{S:def} 
If $A$ is a subset of a ring $R$, we denote respectively by $\ell(A)$ and $\mathrm{r}(A)$ its left annihilator and its right annihilator.
Given a  ring $R$ and a left $R$-module $M$, we say that $M$ is \textbf{P-flat} if, for any $(s,x)\in R\times M$ such that $sx=0$, $x\in \mathrm{r}(s)M$. When $R$ is a domain, $M$ is P-flat if and only if it is torsion-free. As in \cite{Azu87}, we say that $M$ is \textbf{finitely projective} (respectively \textbf{singly projective}) if, for any finitely generated (respectively cyclic) submodule $N$, the inclusion map $N\rightarrow M$ factors through a free module $F$. A finitely projective module is called f-projective in \cite{Jon82}. As in \cite{ZH76} we say that $M$ is \textbf{locally projective} if, for any finitely generated submodule $N$, there exist a free module $F$, an homomorphism $\phi:M\rightarrow F$ and an homomorphism $\pi:F\rightarrow M$ such that $\pi(\phi(x))=x,\ \forall x\in N$. A locally projective module is said to be either a trace module or a universally torsionless module in \cite{Gar76}.
Given a ring $R$, a left $R$-module $M$ and $x\in M$,  the \textbf{content ideal} $\mathrm{c}(x)$ of $x$ in $M$, is the intersection of all right ideals $A$ for which $x\in AM$. We say that $M$ is a \textbf{content module} if $x\in\mathrm{c}(x)M,\ \forall x\in M$.

It is obvious that each locally projective module is finitely projective but the converse doesn't generally hold. For instance, if $R$ is a commutative domain with quotient field $Q\ne R$, then $Q$ is a finitely projective $R$-module: if $N$ is a finitely generated submodule of $Q$, there exists $0\ne s\in R$ such that $sN\subseteq R$, whence the inclusion map $N\rightarrow Q$ factors through $R$ by using the multiplications by $s$ and $s^{-1}$; but $Q$ is not locally projective because the only homormorphism from $Q$ into a free $R$-module is zero.

\begin{proposition} \label{P:pflat}
Let $R$ be a ring. Then:
\begin{enumerate}
\item Each singly
  projective left $R$-module $M$ is P-flat. The converse holds
  if $R$ is a domain.
\item Any P-flat cyclic left module is flat.
\item Each P-flat content left module $M$ is singly projective.
\end{enumerate} 

\end{proposition}
\textbf{Proof.} $(1)$. Let $0\ne x\in M$ and $r\in R$ such that $rx=0$. There
exist a  free module $F$ and two homomorphisms $\phi:Rx\rightarrow
F$ and $\pi:F\rightarrow M$ such that $\pi\circ\phi$ is the inclusion
map $Rx\rightarrow M$. Since $r\phi(x)=0$ and $F$ is free, there
exist $s_1,\dots,s_n\in\mathrm{r}(r)$ and $y_1,\dots,y_n\in F$ such that $\phi(x)=s_1y_1+\dots+s_ny_n$. Then
$x=s_1\pi(y_1)+\dots+s_n\pi(y_n)$. The last assertion is obvious.

$(2)$. Let $C$ be a cyclic left module generated by $x$ and let $A$ be a right ideal. Then each element of $A\otimes_RC$ is of the form $a\otimes x$ for some $a\in A$. If $ax=0$ then $\exists b\in\mathrm{r}(a)$ such that $x=bx$. Therefore $a\otimes x=a\otimes bx=ab\otimes x=0$. Hence $C$ is flat.

$(3)$. Let $x\in M$. Then, since $x\in\mathrm{c}(x)M$ there exist $a_1,\dots,a_n\in\mathrm{c}(x)$ and $x_1,\dots,x_n\in M$ such that $x=a_1x_1+\dots+a_nx_n$. Let $b\in R$ such that $bx=0$. Therefore $x\in\mathrm{r}(b)M$ because $M$ is P-flat. It follows that $\mathrm{c}(x)\subseteq\mathrm{r}(b)$. So, if we put $\phi(rx)=(ra_1,\dots,ra_n)$, then $\phi$ is a well defined  homomorphism which factors the inclusion map $Rx\rightarrow M$ through ${}_RR^n$.
\qed

\begin{theorem}
\label{T:perf} A ring $R$ is left perfect if and only if each flat left module is a content module.
\end{theorem}
\textbf{Proof.} If $R$ is left perfect then each flat left module is projective. Conversely suppose that each flat left module is a content module. Let $(a_k)_{k\in\mathbb{N}}$ be a family of elements of $R$, let $(e_k)_{k\in\mathbb{N}}$ be a basis of a free left module $F$ and let $G$ be the submodule of $F$ generated by $\{e_k-a_ke_{k+1}\mid k\in\mathbb{N}\}$. By \cite[Lemma 1.1]{Bas60} $F/G$ is flat. We put $z_k=e_k+G,\ \forall k\in\mathbb{N}$. Since $F/G$ is content and $z_k=a_kz_{k+1},\ \forall k\in\mathbb{N}$, there exist $c\in R$ and $n\in\mathbb{N}$ such that $z_0=cz_n$ and $\mathrm{c}(z_0)=cR$. It follows that $cR=ca_n\dots a_pR,\ \forall p>n$. Since $z_0=a_0\dots a_{n-1}z_n$, there exists $k>n$ such that $ca_n\dots a_k=a_0\dots a_k$. Consequently $a_0\dots a_kR=a_0\dots a_pR,\ \forall p\geq k$. So, $R$ is left perfect because it satisfies the descending chain condition on principal right ideals by \cite[Theorem P]{Bas60}. \qed

\bigskip
Given a ring $R$ and a left $R$-module $M$, we say that $M$ is \textbf{P-injective} if, for any $(s,x)\in R\times M$ such that $\ell(s)x=0$, $x\in sM$. When $R$ is a domain, $M$ is P-injective if and only if it is divisible. As in \cite{RaRa73}, we say that $M$ is \textbf{finitely injective} (respectively \textbf{FP-injective}) if, for any finitely generated submodule $A$ of a (respectively finitely presented) left module $B$, each homomorphism from $A$ to $M$ extends to $B$. If $M$ is an $R$-module, we put $M^*=\mathrm{Hom}_R(M,R)$.

\begin{proposition} \label{P:ringinj} Let $R$ be a ring. Then: 
\begin{enumerate}
\item If $R$ is a P-injective left module then each  singly
  projective left module is P-injective;
\item If $R$ is a  FP-injective  left module then each  finitely
  projective  left module is FP-injective and a content module;
\item If $R$ is an  injective  module then each  singly
  projective module is  finitely injective and locally projective.
\end{enumerate}
\end{proposition}
\textbf{Proof.} Let $M$ be a left module, $F$ a  free left module and
$\pi:F\rightarrow M$ an epimorphism.

$1.$ Assume that $M$ is singly projective. Let $x\in M$ and $r\in R$
such that $\ell(r)x=0$. There exists a
homomorphism $\phi:Rx\rightarrow F$ such that $\pi\circ\phi$ is the inclusion
map $Rx\rightarrow M$. Since $F$ is P-injective, $\phi(x)=ry$ for some
$y\in F$. Then $x=r\pi(y)$.

$2.$ Assume that $M$ is finitely projective. Let $L$ be a 
finitely generated free left  module, let $N$ be a finitely generated
submodule of $L$ and let $f:N\rightarrow M$ be a homomorphism. Then
$f(N)$ is a finitely generated submodule of $M$. So, there exists a
homomorphism $\phi:f(N)\rightarrow F$ such that $\pi\circ\phi$ is the inclusion
map $f(N)\rightarrow M$. Since $F$ is FP-injective, there exists a
morphism $g:L\rightarrow F$ such that $\phi\circ f$ is the restriction
of $g$ to $N$. Now it is easy to check that $\pi\circ g$ is the
restriction of $f$ to $N$.

Let $x\in M$. There exists a
homomorphism $\phi:Rx\rightarrow F$ such that $\pi\circ\phi$ is the inclusion
map $Rx\rightarrow M$. Let $\{e_i\mid i\in I\}$ be a basis of $F$. There exist a finite subset $J$ of $I$ and a family $(a_i)_{i\in J}$ of elements of $R$ such that $\phi(x)=\sum_{i\in J}a_ie_i$. 	Let $A$ be the right ideal generated by $(a_i)_{i\in J}$. Then $(0:x)=(0:\phi(x))=\ell(A)$. Let $B$ be a right ideal such that $x\in BM$. Then $x=\sum_{k=1}^pb_kx_k$ where $b_k\in B$ and $x_k\in M$, $\forall k,\ 1\leq k\leq p$. Let $N$ be the submodule of $M$ generated by $\{\pi(e_i)\mid i\in J\}\cup\{x_k\mid 1\leq k\leq p\}$. Thus there exists a
homomorphism $\varphi:N\rightarrow F$ such that $\pi\circ\varphi$ is the inclusion
map $N\rightarrow M$. Therefore there exist a finite subset $K$ of $I$ and two families $\{d_{k,j}\mid 1\leq k\leq p,\ j\in K\}$ and $\{c_{i,j}\mid(i,j)\in J\times K\}$ of elements of $R$ such that $\varphi(\pi(e_i))=\sum_{j\in K}c_{i,j}e_j$, $\forall i\in J$ and $\varphi(x_k)=\sum_{j\in K}d_{k,j}e_j$, $\forall k,\ 1\leq k\leq p$. It follows that $\varphi(x)=\sum_{j\in K}(\sum_{i\in J}a_ic_{i,j})e_j=\sum_{j\in K}(\sum_{k=1}^pb_kd_{k,j})e_j$. So, $\sum_{i\in J}a_ic_{i,j}=\sum_{k=1}^pb_kd_{k,j},\ \forall j\in K$. Let $A'$ be the right ideal generated by $\{\sum_{i\in J}a_ic_{i,j}\mid j\in K\}$. Then $A'\subseteq A$ and $A'\subseteq B$. Moreover, $\ell(A)=(0:x)=(0:\varphi(x))=\ell(A')$. By \cite[Corollary 2.5]{Ja73} $A=A'$. So, $A\subseteq B$. We conclude that $\mathrm{c}(x)=A$ and $M$ is a content module.

$3.$ Let $M$ be a singly projective module and $x\in M$.  So, there exists a
homomorphism $\phi:Rx\rightarrow F$ such that $\pi\circ\phi$ is the inclusion
map $Rx\rightarrow M$. Since $F$ is finitely injective, we can extend
$\phi$ to $M$. By using a basis of $F$ we deduce that $x=\sum_{k=1}^n\phi_k(x)x_k$ where
$\phi_k\in M^*$ and $x_k\in M,\ \forall k,\ 1\leq k\leq n$. Hence $M$ is
locally projective by \cite[Theorem 3.2]{Gar76} or \cite[Theorem 2.4]{ZH76}.  By a similar proof as in $(2)$, we show that $M$ is finitely injective, except that $L$ is not necessarily a finitely
generated free module. \qed

\bigskip
A short exact sequence of left $R$-modules $0\rightarrow N\rightarrow M\rightarrow L\rightarrow 0$ is \textbf{pure} if it remains exact when tensoring it with any right $R$-module. We say that $N$ is a pure submodule of $M$. This property holds if $L$ is flat.
\begin{lemma} \label{L:free}
Let $R$ be a local ring, let $P$ be its maximal ideal and let $N$ be a flat left $R$-module. Assume that
$N$ is generated by a family $(x_i)_{i\in I}$ of elements of $N$ such
that $(x_i+PN)_{i\in I}$ is a basis of $N/PN$. Then $N$ is free.
\end{lemma}
\textbf{Proof.} Let $(e_i)_{i\in I}$ be a basis of a
free left module $F$, let $\alpha:F\rightarrow N$ be the homomorphism
defined by $\alpha(e_i)=x_i,\ \forall i\in I$ and let $L$ be the kernel of
$\alpha$. It is easy to check that $L\subseteq PF$. Let  $y\in L$
. We have $y=\sum_{i\in J}a_ie_i$
where  $J$ is a finite subset of $I$ and $a_i\in P,\ \forall
i\in J$. Since $L$ is a pure submodule of $F$, $\forall i\in J$ there exists
$y_i\in L$ such that $\sum_{i\in J}a_ie_i=\sum_{i\in J}a_iy_i$. We
have $y_i=\sum_{j\in J_i}b_{i,j}e_j$ where $J_i$ is a finite subset of
  $I$, $b_{i,j}\in P,\ \forall (i,j)\in J\times J_i$. Let $K=J\cup (\cup_{i\in J}J_i)$. If $i\in K\setminus J$ we put $a_i=0$ and $a_{i,j}=0,\ \forall j\in K$, and if $j\in K\setminus J_i$ we put $a_{i,j}=0$ too.
  We get $\sum_{ i\in K}a_ie_i=\sum_{j\in
   K}(\sum_{i\in K}a_ib_{i,j})e_j$. It follows that $a_j=\sum_{i\in
   K}a_ib_{i,j}$. So, if $A$ is the right ideal generated by $\{a_i\mid i\in
   K\}$, then $A=AP$. By Nakayama lemma $A=0$, whence $F\cong N$. \qed

\bigskip
A left $R$-module is said to be a \textbf{Mittag-Leffler} module if, for each index set $\Lambda$, the natural homomorphism $R^{\Lambda}\otimes_RM\rightarrow M^{\Lambda}$ is injective. The following lemma is a slight generalization of \cite[Proposition 2.3]{CoPe70}.

\begin{lemma} \label{L:flat-fproj} Let $R$ be a subring of a ring $S$ and let $M$ be a flat left $R$-module. Assume that $S\otimes_RM$ is finitely projective over $S$. Then $M$ is finitely projective.
\end{lemma}
\textbf{Proof.} By \cite[Proposition 2.7]{Jon82} a module is finitely projective if and only if it is a flat Mittag-Leffler module. So we do as in the proof of \cite[Proposition 2.3]{CoPe70}. \qed

\bigskip
From this lemma and \cite[Proposition 2.7]{Jon82} we deduce the following proposition. We can also see .
\begin{proposition} \label{P:Flat=F-proj} Let $R$ be a subring of a left perfect ring $S$. Then each flat left $R$-module is finitely projective.
\end{proposition}

\begin{proposition} \label{P:locali} 
Let $R$ be a commutative ring and let $S$ be a multiplicative subset of $R$. Then:
\begin{enumerate}
\item For each singly (respectively finitely, locally) projective
  $R$-module $M$, $S^{-1}M$ is singly (respectively finitely, locally) projective over
  $S^{-1}R$;
\item Let $M$ be
  a singly (respectively finitely) projective $S^{-1}R$-module. If $S$ contains no
  zero-divisors then $M$ is singly (respectively finitely) projective
  over $R$.
\end{enumerate}
\end{proposition}
\textbf{Proof.} $(1)$. We assume that $M\ne 0$. Let $N$ be a cyclic (respectively finitely generated)
submodule of $S^{-1}M$. Then there exists a cyclic (respectively finitely generated)
submodule $N'$ of $M$ such that
$S^{-1}N'=N$. There exists a free $R$-module $F$, a morphism
$\phi:N'\rightarrow F$ and a morphism $\pi:F\rightarrow M$ such that
$(\pi\circ\phi)(x)=x$ for each $x\in N'$. It follows that
$(S^{-1}\pi\circ S^{-1}\phi)(x)=x$ for each $x\in N$. We get that
$S^{-1}M$ is singly (respectively finitely) projective over $R$. We do a similar proof to show that $S^{-1}M$ is locally projective if $M$ is locally projective.

$(2)$ By Lemma~\ref{L:flat-fproj} $M$ is finitely projective over $R$ if it is finitely projective over $S^{-1}R$. It is easy to check that $M$ is singly projective over $R$ if it is singly projective over $S^{-1}R$.
\qed

\bigskip
If $R$ is a subring of a ring $Q$ which is either left perfect or left noetherian, then then each flat left $R$-module is finitely projective by \cite[Corollary 7]{She91}. We don't know if the converse holds. However we have the following results:
\begin{theorem} \label{T:Q-FP-inj} Let $R$ be a commutative ring with a self FP-injective quotient ring $Q$. Then each flat $R$-module is finitely projective if and only if $Q$ is perfect.
\end{theorem}
\textbf{Proof.} "Only if" requires a proof. Let $M$ be a flat $Q$-module. Then $M$ is flat over $R$ and it follows that $M$ is finitely projective over $R$. By Proposition~\ref{P:locali}(1) $M\cong Q\otimes_RM$ is finitely projective over $Q$. From Proposition~\ref{P:ringinj} we deduce that each flat $Q$-module is content. We conclude by Theorem~\ref{T:perf} \qed

\begin{theorem} \label{T:Flat=f-pro} Let $R$ be a commutative ring with a Von Neumann regular quotient ring $Q$. Then the following conditions are equivalent:
\begin{enumerate}
\item $Q$ is semi-simple;
\item each flat $R$-module is finitely projective;
\item each flat $R$-module is singly projective.
\end{enumerate}   
\end{theorem}
\textbf{Proof.} $(1)\Rightarrow (2)$ is an immediate consequence of \cite[Corollary 7]{She91} and $(2)\Rightarrow (3)$ is obvious. 

$(3)\Rightarrow (1)$. First we show that each $Q$-module $M$ is singly projective. Every $Q$-module $M$ is flat over $Q$ and $R$. So, $M$ is singly projective over $R$. It follows that $M\cong Q\otimes_RM$ is singly projective over $Q$ by Proposition~\ref{P:locali}(1). Now let $A$ be an ideal of $Q$. Since $Q/A$ is singly projective, it is projective. So, $Q/A$ is finitely presented over $Q$ and $A$ is a finitely generated ideal of $Q$. Hence $Q$ is semi-simple. 
\qed

\section{$\pi$-coherence and PP-rings}
\label{S:baer}

As in \cite{ZH76} we say  that a ring $R$ is left \textbf{strongly coherent} if each product of
 locally projective right modules is locally projective and as in \cite{Cam90}  $R$ is said to be right \textbf{$\pi$-coherent} if, for each index set $\Lambda$, every finitely generated submodule of $R_R^{\Lambda}$ is finitely presented.
 
\begin{theorem}
\label{T:picoh} Let $R$ be a commutative ring. Then the following conditions are equivalent:
\begin{enumerate}
\item $R$ is $\pi$-coherent;
\item for each index set $\Lambda$, $R^{\Lambda}$ is finitely projective;
\item each product of finitely projective modules is finitely projective.
\end{enumerate} 
\end{theorem}
\textbf{Proof.}
$(1)\Rightarrow (2)$. Let $N$ be a finitely generated submodule of $R^{\Lambda}$. There exist a free module $F$ and an epimorphism $\pi$ from $F$ into $R^{\Lambda}$. It is obvious that $R$ is coherent. Consequently $R^{\Lambda}$ is flat. So $\ker\ \pi$ is a pure submodule of $F$. Since $N$ is finitely presented it follows that there exists $\phi:N\rightarrow F$ such that $\pi\circ\phi$ is the inclusion map from $N$ into $R^{\Lambda}$.

$(2)\Rightarrow (1)$. Since $R^{\Lambda}$ is flat for each index set $\Lambda$, $R$ is coherent. Let $\Lambda$ be an index set and let $N$ be a finitely generated submodule of $R^{\Lambda}$. The finite projectivity of $R^{\Lambda}$ implies that $N$ is isomorphic to a submodule of a free module of finite rank. Hence $N$ is finitely presented.

It is obvious that $(3)\Rightarrow (2)$.

$(2)\Rightarrow (3)$. Let $\Lambda$ be an index set, let $(M_{\lambda})_{\lambda\in\Lambda}$ be a family of finitely projective modules and let $N$ be a finitely generated submodule of $M=\prod_{\lambda\in\Lambda}M_{\lambda}$. For each $\lambda\in\Lambda$, let $N_{\lambda}$ be the image of $N$ by the canonical map $M\rightarrow M_{\lambda}$. We put $N'=\prod_{\lambda\in\Lambda}N_{\lambda}$. So, $N\subseteq N'\subseteq M$. For each $\lambda\in\Lambda$ there exists a free module $F_{\lambda}$ of finite rank such that the inclusion map $N_{\lambda}\rightarrow M_{\lambda}$ factors through $F_{\lambda}$. It follows that the inclusion map $N\rightarrow M$ factors through $\prod_{\lambda\in\Lambda}F_{\lambda}$ which is isomorphic to $R^{\Lambda'}$ for some index set $\Lambda'$. Now the monomorphism $N\rightarrow R^{\Lambda'}$ factors through a free module $F$. It is easy to conclude that the inclusion map $N\rightarrow M$ factors through $F$ and that $M$ is finitely projective. \qed

\bigskip 
By using \cite[Theorem 4.2]{ZH76} and Proposition~\ref{P:ringinj}, we deduce the following corollary:
\begin{corollary}
\label{C:LocFinPro} Every strongly coherent commutative ring $R$ is $\pi$-coherent and the converse holds if $R$ is self injective.
\end{corollary}
\begin{proposition}
\label{P:pilocali} Let $R$ be a $\pi$-coherent commutative ring and let $S$ be a multiplicative subset of $R$. 	Assume that $S$ contains no zero-divisors. Then $S^{-1}R$ is $\pi$-coherent.
\end{proposition}
\textbf{Proof.} Let $M$ be a finitely generated $S^{-1}R$-module. By \cite[Theorem 1]{Cam90} we must prove that $\mathrm{Hom}_{S^{-1}R}(M,{S^{-1}R})$ is finitely generated on $S^{-1}R$. There exists a finitely generated $R$-submodule $N$ of $M$ such that $S^{-1}N\cong M$. The following sequence
\[0\rightarrow N^*\rightarrow\mathrm{Hom}_R(N,S^{-1}R)\rightarrow\mathrm{Hom}_R(N,S^{-1}R/R)\]
is exact. Since $N$ is finitely generated and $S^{-1}R/R$ is $S$-torsion, $\mathrm{Hom}_R(N,S^{-1}R/R)$ is $S$-torsion too. So, $\mathrm{Hom}_{S^{-1}R}(M,S^{-1}R)\cong\mathrm{Hom}_R(N,S^{-1}R)\cong S^{-1}N^*$. By \cite[Theorem 1]{Cam90} $N^*$ is finitely generated. Hence $\mathrm{Hom}_{S^{-1}R}(M,S^{-1}R)$ is finitely generated over $S^{-1}R$. \qed

\begin{theorem}
\label{T:semihered} Let $R$ be a commutative semihereditary ring and let $Q$ be its quotient ring. Then the following conditions are equivalent:
\begin{enumerate}
\item $R$ is $\pi$-coherent;
\item $Q$ is self injective;
\end{enumerate}
Moreover, when these conditions are satisfied, each singly projective $R$-module is finitely projective.
\end{theorem}
\textbf{Proof.} $(1)\Rightarrow (2)$. By Proposition~\ref{P:pilocali} $Q$ is $\pi$-coherent. We know that $Q$ is Von Neumann regular. It follows from \cite[Theorem 2]{Kob84} that $Q$ is self injective.

$(2)\Rightarrow (1)$. Let $(M_i)_{i\in I}$ be a family of finitely projective $R$-modules, where $I$ is an index set, and let $N$ be a finitely generated submodule of $\prod_{i\in I}M_i$. Then $N$ is flat. Since $N$ is a submodule of $\prod_{i\in I}Q\otimes_RM_i$, $Q\otimes_RN$ is isomorphic to a finitely generated $Q$-submodule of $\prod_{i\in I}Q\otimes_RM_i$. It follows that $Q\otimes_RN$ is a projective $Q$-module. Hence $N$ is projective by
\cite[Proposition 2.3]{CoPe70}. We conclude by Theorem~\ref{T:picoh}.

Let $M$ be a singly projective $R$-module and let $N$ be a finitely generated submodule of $M$. Then $Q\otimes_RM$ is finitely projective over $Q$ by Propositions \ref{P:locali}(1) and \ref{P:ringinj}. It follows that $Q\otimes_RN$ is projective over $Q$. Hence $N$ is projective by
\cite[Proposition 2.3]{CoPe70}. \qed

\begin{proposition}
\label{P:cont=sproj} Let $R$ be a Von Neumann regular ring. Then a right $R$-module is  content if and only if it is singly projective.
\end{proposition}
\textbf{Proof.} By Proposition~\ref{P:pflat}$(3)$ it remains to show that each singly projective right module $M$ is content. Let $m\in M$. Then $mR$ is projective because it is isomorphic to a finitely generated submodule of a free module. So, $mR$ is content. For each left ideal $A$, $mR\cap MA=mA$ because $mR$ is a pure submodule of $M$. Hence $M$ is content. \qed

\bigskip
A topological space $X$ is said to be \textbf{extremally disconnected} if every open set has an open closure. Let $R$ be a ring. We say that $R$ is a right \textbf{Baer ring} if for any subset $A$ of $R$, $\mathrm{r}(A)$ is generated by an idempotent. The ring $R$ defined in  \cite[Example 4.4 ]{ZH76} is not self injective and satisfies the conditions of the following theorem. 

\begin{theorem}
\label{T:regring} Let $R$ be a Von Neumann regular ring. Then the following conditions are equivalent:
\begin{enumerate}
\item Each product of singly projective right modules is singly projective;
\item Each product of content right modules is content;
\item $R_R^R$ is singly projective;
\item $R_R^R$ is a content module;
\item $R$ is a right Baer ring;
\item The intersection of each family of finitely generated left ideals is finitely generated too;
\item For each cyclic left module $C$, $C^*$ is finitely generated.
\end{enumerate}
Moreover, when $R$ is commutative, these conditions are equivalent to the following: $\mathrm{Spec}\ R$ is extremally disconnected.

\end{theorem}
\textbf{Proof.} The conditions $(2),(4),(6)$ are equivalent by \cite[Theorem 5.15]{Gar76}. By Proposition~\ref{P:cont=sproj} $(4)\Leftrightarrow (3)$ and $(1)\Leftrightarrow (2)$. It is easy to check that $(5)\Leftrightarrow (7)$.

$(3)\Rightarrow (5)$. Let $A\subseteq R$ and let $x=(a)_{a\in A}\in R_R^A$. So, $\mathrm{r}(A)=(0:x)$. Then $xR$ is projective because it is isomorphic to a submodule of a free module. Thus $\mathrm{r}(A)=eR$, where $e$ is an idempotent.

$(5)\Rightarrow (1)$. Let $(M_i)_{i\in I}$ be a family of singly projective right modules and $m=(m_i)_{i\in I}$ be an element of $M=\prod_{i\in I}M_i$. For each $i\in I$, there exists an idempotent $e_i$ such that $(0:m_i)=e_iR$. Let $e$ be the idempotent which satisfies $eR=\mathrm{r}(\{1-e_i\mid i\in I\})$. Then $eR=(0:m)$, whence $mR$ is projective.

If $R$ is commutative and reduced, then the closure of $D(A)$, where $A$ is an ideal of $R$, is $V((0:A))$. So, $\mathrm{Spec}\ R$ is extremally disconnected if and only if, for each ideal $A$ there exists an idempotent $e$ such that $V((0:A))=V(e)$. This last equality holds if and only if $(0:A)=Re$ because $(0:A)$ and $Re$ are semiprime since $R$ is reduced. Consequently $\mathrm{Spec}\ R$ is extremally disconnected if and only if $R$ is Baer. The proof is now complete.
\qed

\bigskip
Let $R$ be a ring. We say that $R$ is a right \textbf{PP-ring} if  any principal right ideal is projective.
 
\begin{lemma} \label{L:cyclicPP}
Let $R$ be a right PP-ring. Then each cyclic submodule of a free right module is projective.
\end{lemma}
\textbf{Proof.} Let $C$ be a cyclic submodule of a free right module $F$. We may assume that $F$ is finitely generated 
by the basis $\{e_1,\dots,e_n\}$. Let $p:F\rightarrow R$ be the homomorphism defined by $p(e_1r_1+\dots+e_nr_n)=r_n$ where $r_1,\dots,r_n\in R$. Then $p(C)$ is a principal right ideal. Since $p(C)$ is projective, $C\cong C'\oplus p(C)$ where $C'=C\cap\mathrm{ker}\ p$. So $C'$ is a cyclic submodule of the free right module generated by $\{e_1,\dots,e_{n-1}\}$. We complete the proof by induction on $n$. \qed

\begin{theorem} \label{T:baer}
Let $R$ be a commutative PP-ring and let $Q$ be its quotient ring. Then the following conditions are equivalent:
\begin{enumerate}
\item Each product of singly projective modules is singly projective;
\item $R^R$ is singly projective;
\item $R$ is a  Baer ring;
\item $Q$ satisfies the equivalent conditions of Theorem~\ref{T:regring};
\item For each cyclic module $C$, $C^*$ is finitely generated;
\item $\mathrm{Spec}\ R$ is extremally disconnected;
\item $\mathrm{Min}\ R$ is extremally disconnected.
\end{enumerate}
\end{theorem}
\textbf{Proof.} It is obvious that $(1)\Rightarrow (2)$. It is easy to check that $(3)\Leftrightarrow (5)$. We show that $(2)\Rightarrow (3)$ as we proved $(3)\Rightarrow (5)$ in Theorem~\ref{T:regring}, by using Lemma~\ref{L:cyclicPP}.

$(5)\Rightarrow (4)$. Let $C$ be a cyclic $Q$-module. We do as in proof of Proposition~\ref{P:pilocali} to show that $\mathrm{Hom}_Q(C,Q)$ is finitely generated over $Q$.

$(4)\Rightarrow (1)$. Let $(M_i)_{i\in I}$ be a family of singly projective right modules and let $N$ be a cyclic submodule of $M=\prod_{i\in I}M_i$. Since $R$ is PP, $N$ is a P-flat module. By Proposition~\ref{P:pflat} $N$ is flat.
We do as in the proof of $(2)\Rightarrow (1)$ of Theorem~\ref{T:semihered} to show that $N$ is projective.
 
$(3)\Leftrightarrow (6)$ is shown in the proof of Theorem~\ref{T:regring}.
 
$(4)\Leftrightarrow (7)$ holds because $\mathrm{Spec}\ Q$ is homeomorphic to $\mathrm{Min}\ R$. \qed


\section{Flat modules}
\label{S:flat} 
Let $M$ be a non-zero module over a commutative ring $R$. As in
\cite[p.338]{FuSa01} we set:
\[M_{\sharp}=\{s\in R\mid \exists 0\ne x\in M\ \mathrm{such\ that}\
sx=0\}\quad\mathrm{and}\quad M^{\sharp}=\{s\in R\mid sM\subset M\}.\] 
Then $R\setminus M_{\sharp}$ and $R\setminus M^{\sharp}$ are multiplicative subsets of $R$.

\begin{lemma} \label{L:prime}
Let $M$ be a non-zero P-flat $R$-module over a commutative ring $R$. Then $M_{\sharp}\subseteq R_{\sharp}\cap
M^{\sharp}.$
\end{lemma}
\textbf{Proof.} Let $0\ne s\in M_{\sharp}$. Then there exists $0\ne
x\in M$ such that $sx=0$. Since $M$ is P-flat, we have $x\in
(0:s)M$. Hence $(0:s)\ne 0$ and $s\in R_{\sharp}$.

Suppose that $M_{\sharp}\nsubseteq M^{\sharp}$ and let $s\in
M_{\sharp}\setminus M^{\sharp}$. Then $\exists 0\ne x\in M$ such that
$sx=0$. It follows that $x=t_1y_1+\dots+t_py_p$ for some $y_1,\dots,y_p\in M$ and $t_1,\dots,t_p\in
(0:s)$. Since $s\notin M^{\sharp}$ we have $M=sM$. So $y_k=sz_k$ for some
$z_k\in M$, $\forall k,\ 1\leq k\leq p$. We get $x=t_1sz_1+\dots+t_psz_p=0$. Whence a contradiction. \qed

\bigskip
Now we assume that $R$ is a commutative ring.
 An $R$-module $M$ is said to be \textbf{uniserial} if its set of submodules is totally ordered by inclusion and  $R$ is a \textbf{valuation ring} if it is uniserial as $R$-module. If $M$ is a module over a valuation ring $R$ then $M_{\sharp}$ and $M^{\sharp}$ are prime ideals of $R$.
In the sequel, if $R$ is a valuation ring, we denote by $P$ its maximal ideal and we put $Z=R_{\sharp}$ and $Q=R_Z$. Since each finitely generated ideal of a valuation ring $R$ is principal, it follows that any P-flat $R$-module is flat.
\begin{proposition} \label{P:hullflat}
Let $R$ be a valuation ring, let $M$ be a flat $R$-module and let $E$ be its injective hull. Then $E$
is flat.
\end{proposition}
\textbf{Proof.} Let $x\in E\setminus M$ and $r\in R$ such that
$rx=0$. There exists $a\in R$ such that $0\ne ax\in M$. From $ax\ne 0$
and $rx=0$ we deduce that $r=ac$ for some $c\in R$. Since $cax=0$ and
$M$ is flat we have $ax=by$ for some $y\in M$ and $b\in (0:c)$. From
$bc=0$ and $ac=r\ne 0$ we get $b=at$ for some $t\in R$. We have
$a(x-ty)=0$. Since $at=b\ne 0$, $(0:t)\subset Ra$. So $(0:t)\subseteq
(0:x-ty)$. The injectivity of $E$ implies that there exists $z\in E$
such that $x=t(y+z)$. On the other hand $tr=tac=bc=0$, so $t\in (0:r)$. 
\qed

\bigskip

In the sequel, if $J$ is a prime ideal of $R$ we denote by $0_J$ the
kernel of the natural map: $R\rightarrow R_J$.

\begin{proposition} \label{P:ann}
Let $R$ be a valuation ring and let $M$ be a non-zero flat $R$-module. Then:
\begin{enumerate}
\item If $M_{\sharp}\subset Z$ we have $\mathrm{ann}(M)=0_{M_{
      \sharp}}$ and $M$ is an $R_{M_{\sharp}}$-module;
\item If $M_{\sharp}=Z$, $\mathrm{ann}(M)=0$ if $M_Z\ne ZM_Z$ and $\mathrm{ann}(M)=(0:Z)$ if
      $M_Z=ZM_Z$. In this last case, 
      $M$ is a $Q$-module.
\end{enumerate}
\end{proposition}
\textbf{Proof.} Observe that the natural map $M\rightarrow
M_{M_{\sharp}}$ is a monomorphism. First we assume that $R$ is
self FP-injective and $P=M_{\sharp}$. So $M^{\sharp}=P$ by Lemma~\ref{L:prime}. If $M\ne PM$ let $x\in
M\setminus PM$. Then $(0:x)=0$ else $\exists r\in R,\ r\ne 0$ such that $x\in (0:r)M\subseteq PM$. If $M=PM$ then $P$ is not
finitely generated else $M=pM$, where $P=pR$, and $p\notin M^
{\sharp}=P$. If $P$ is not faithful then
$(0:P)\subseteq\mathrm{ann}(M)$. Thus $M$ is flat over $R/(0:P)$. So
we can replace $R$ with $R/(0:P)$ and assume that $P$ is
faithful. Suppose $\exists 0\ne r\in P$ such that $rM=0$. Then
$M=(0:r)M$. Since $(0:r)\ne P$, let $t\in P\setminus (0:r)$. Thus
$M=tM$ and $t\notin M^{\sharp}=P$. Whence a contradiction. So $M$ is
faithful or $\mathrm{ann}(M)=(0:P)$.

Return to the general case. We put $J=M_{\sharp}$. 

If $J\subset Z$
then $R_J$ is coherent and self FP-injective by \cite[Theorem 11]{Couch03}. In this case  $JR_J$ is principal or
faithful. So $M_J$ is faithful over $R_J$,
whence $\mathrm{ann}(M)=0_J$. Let $s\in R\setminus J$. There exists $t\in
Zs\setminus J$. It is easy to check that $\forall a\in R,\ (0:a)$ is also an ideal of $Q$. On the other hand, $\forall a\in Q,\ Qa=(0:(0:a))$ because $Q$ is self FP-injective. It follows that $(0:s)\subset (0:t)$. Let $r\in (0:t)\setminus (0:s)$. Then $r\in 0_J$. So,  $rM=0$. Hence
$M=(0:r)M=sM$. Therefore the
multiplication by $s$ in $M$ is bijective for each $s\in R\setminus J$. 

Now suppose that $J=Z$. Since $Q$ is self FP-injective then $M$ is faithful or $\mathrm{ann}(M)=(0:Z)$. Let $s\in R\setminus Z$. Thus $Z\subset Rs$ and $sZ=Z$. It follows that $ZM_Z=ZM$. So, $M$ is a $Q$-module if $ZM_Z=M_Z$. \qed

\bigskip
 When $R$ is a valuation ring, $N$ is a pure submodule of $M$ if $rN=rM\cap N,\ \forall r\in R$. 
\begin{proposition} \label{P:unis}
 Let $R$ be a valuation ring and let $M$ be a non-zero flat $R$-module such that
  $M_{M_{\sharp}}\ne M_{\sharp}M_{M_{\sharp}}$. Then $M$ contains a
  non-zero pure uniserial submodule.
\end{proposition}
\textbf{Proof.} Let $J=M_{\sharp}$ and $x\in M_J\setminus JM_J$. If
$J\subset Z$ then $M$ is a module over $R/0_J$ and $J/0_J$ is the
subset of zero-divisors of $R/0_J$. So, after replacing $R$ with
$R/0_J$ we may assume that $Z=J$. If
$rx=0$ then $x\in (0:r)M_Z\subseteq ZM_Z$ if $r\ne 0$. Hence $Qx$ is
faithful over $Q$ which is FP-injective. So $V=Qx$ is a pure
submodule of $M_Z$. We put $U=M\cap V$. Thus $U$ is uniserial and
$U_Z=V$. Then $M/U$ is a submodule of $M_Z/V$, and this last module is
flat. Let $z\in M/U$ and $0\ne r\in R$ such that $rz=0$. Then
$z=as^{-1}y$ where $s\notin Z$, $a\in (0:r)\subseteq Z$ and $y\in
M/U$. It follows that $a=bs$ for some $b\in R$ and $sbr=0$. So $b\in (0:r)$ and
$z=by$. Since $M/U$ is flat, $U$ is a pure submodule of $M$. \qed

\begin{proposition} \label{P:freepure}
Let $R$ be a valuation ring and let $M$ be a flat $R$-module. Then $M$ contains a pure free submodule $N$ such that $M/PM\cong N/PN$.
\end{proposition}
\textbf{Proof.}
 Let $(x_i)_{i\in
  I}$ be a family of elements of $M$ such that $(x_i+PM)_{i\in I}$ is
a basis of $M/PM$ over $R/P$, and let $N$ be the submodule of $M$
generated by this family. If we show that $N$ is a pure
submodule of $M$, we deduce that $N$ is flat. It follows that $N$ is free by Lemma~\ref{L:free}. Let $x\in M$ and $r\in
R$ such that $rx\in N$. Then $rx=\sum_{i\in J}a_ix_i$ where $J$ is a
finite subset of $I$ and $a_i\in R,\ \forall i\in J$. Let $a\in R$
such that $Ra=\sum_{i\in J}Ra_i$. It follows that, $\forall i\in J$, there exists $u_i\in
R$ such that $a_i=au_i$ and there is at least one $i\in
J$ such that $u_i$ is a unit. Suppose that $a\notin Rr$. Thus there
exists $c\in P$ such that $r=ac$. We get that $a(\sum_{i\in
  J}u_ix_i-cx)=0$. Since $M$ is flat, we deduce that $\sum_{i\in
  J}u_ix_i\in PM$. This contradicts that $(x_i+PM)_{i\in I}$ is
a basis of $M/PM$ over $R/P$. So, $a\in Rr$. Hence $N$ is a pure submodule. \qed

\section{Singly projective modules}
\label{S:sinpro}

\begin{lemma} \label{L:sinval} Let $R$ be a valuation ring. Then
a non-zero $R$-module $M$ is singly projective if and only if for
  each $x\in M$ there exists $y\in M$ such that $(0:y)=0$ and $x\in
  Ry$. Moreover $M_{\sharp}=Z$ and $M_Z\ne ZM_Z$.
\end{lemma}
\textbf{Proof.}  Assume that $M$ is singly projective and let $x\in
M$. There exist a free module $F$, a morphism $\phi:Rx\rightarrow F$
and a morphism $\pi:F\rightarrow M$ such that
$(\pi\circ\phi)(x)=x$. Let $(e_i)_{i\in I}$ be a basis of $F$. Then
$\phi(x)=\sum_{i\in J}a_ie_i$ where $J$ is a finite subset of $I$ and
$a_i\in R,\ \forall i\in J$. There exists $a\in R$ such that
$\sum_{i\in J}Ra_i=Ra$. Thus, $\forall i\in J$ there exists $u_i\in R$
such $a_i=au_i$. We put $z=\sum_{i\in J}u_ie_i$. Then
$\phi(x)=az$. Since there is at least one index $i\in J$ such that $u_i$
is a unit, then $(0:z)=0$. It follows that $(0:\phi(x))=(0:a)$. But
$(0:x)=(0:\phi(x))$ because $\phi$ is a monomorphism. We have
$x=a\pi(z)$. So, by \cite[Lemma 2]{Couch03} $(0:\pi(z))=a(0:x)=a(0:a)=0$. The converse and the last assertion are
obvious.
\qed

\bigskip


Let $R$ be a valuation ring and let $M$ be a non-zero
$R$-module. A submodule $N$ of $M$ is said to be \textbf{pure-essential} if it is a pure submodule and if $0$ is the only submodule $K$ satisfying $N\cap K=0$ and $(N+K)/K$ is a pure submodule of $M/K$. An $R$-module $E$ is said to be \textbf{pure-injective} if for any pure-exact sequence $0\rightarrow N\rightarrow M\rightarrow L\rightarrow 0$, the following sequence is exact:
\[0\rightarrow\mathrm{Hom}_R(L,E)\rightarrow\mathrm{Hom}_R(M,E)\rightarrow\mathrm{Hom}_R(N,E)\rightarrow\ 0.\]
We say that $E$ is a \textbf{pure-injective hull} of $K$ if $E$ is pure-injective and $K$ is a pure-essential submodule of $E$. We say that $R$ is \textbf{maximal} if every family of cosets $\{a_i+L_i\mid i\in I\}$ with the finite intersection property has a non-empty intersection (here $a_i\in R,\ L_i$ denote ideals of $R$, and $I$ is an arbitrary index set).

\begin{proposition} \label{P:content}

Let $R$ be a valuation ring and let $M$ be a non-zero
$R$-module. Then the following conditions are equivalent:
\begin{enumerate}
\item $M$ is a flat content module;
\item $M$ is flat and  contains a pure-essential free submodule.
\end{enumerate}
Moreover, if these conditions are satisfied, then any element of
$M$ is contained in a pure cyclic free submodule $L$ of
$M$. If $R$ is maximal then $M$ is locally projective.

\end{proposition}
\textbf{Proof.} 
$(1)\Rightarrow (2)$. Let $0\ne x\in M$. Then $x=\sum_{1\leq i\leq n}a_ix_i$ where $a_i\in\mathrm{c}(x)$ and $x_i\in M,\ \forall i,\ 1\leq i\leq n$. Since $R$ is a valuation ring $\exists a\in R$ such that $Ra=Ra_1+\dots+Ra_n$.  So, we get that $\mathrm{c}(x)=Ra$ and $x=ay$ for some $y\in M$. Thus $y\notin PM$ else $\mathrm{c}(x)\subset Ra$. So $PM\ne M$ and we can apply Proposition~\ref{P:freepure}.
 It remains to show that $N$ is a pure-essential submodule of
$M$. Let $x\in M$ such that $Rx\cap N=0$ and $N$ is a pure submodule
of $M/Rx$. There exist $b\in R$ and $y\in M\setminus PM$ such that
$x=by$. Since $M=N+PM$, we have $y=n+pm$ where $n\in N,\ m\in M$ and
$p\in P$. Then $n\notin PN$ and $bpm=-bn+x$. Since $N$ is pure in $M/Rx$
there exist $n'\in N$ and $t\in R$ such that $bpn'=-bn+tx$. We get
that $b(n+pn')\in N\cap Rx=0$. So $b=0$ because $n+pn'\notin PN$. Hence
$x=0$. 

$(2)\Rightarrow (1)$. First we show that $M$ is a content module if each element $x$ of $M$ is
of the form $s(y+cz)$, where $s\in R,\ y\in N\setminus PN,\ c\in P$ and $z\in
M$. Since $N$ is a pure submodule, $PM\cap N=PN$ whence $y\notin PM$. If $x=stw$ with $t\in P$ and $w\in M$ we get that $s(y+cz-tw)=0$ whence $y\in PM$ because $M$ is flat. This is a contradiction. Consequently $\mathrm{c}(x)=Rs$ and $M$ is content. Now we prove that each element $x$ of $M$ is
of the form $s(y+cz)$, where $s\in R,\ y\in N\setminus PN,\ c\in P$ and $z\in
M$. If $x\in N$, then we check this property by using a basis of
$N$. Suppose $x\notin N$ and $Rx\cap N\ne 0$. There exists $a\in P$
such that $0\ne ax\in N$. Since $N$ is pure, there exists $y'\in N$
such  that $ax=ay'$. We get $x=y'+bz$ for some $b\in (0:a)$ and $z\in
M$, because $M$ is flat. We have $y'=sy$ with $s\in R$ and $y\in
N\setminus PN$. Since $as\ne 0,\ b=sc$ for some $c\in P$. Hence
$x=s(y+cz)$. Now suppose that $Rx\cap N=0$. Since $N$ is
pure-essential in $M$,  there exist $r\in R$ and $m\in M$ such that
$rm\in N+Rx$ and $rm\notin rN+Rx$. Hence $rm=n+tx$ where $n\in N$ and $
t\in R$. Thus $n=by'$ where $b\in R$ and $y'\in N\setminus PN$. Then
$b\notin rR$. So, $r=bc$ for some $c\in P$. We get $bcm=by'+tx$. If
$t=bd$ for some $d\in P$ then $b(cm-y'-dx)=0$. Since $M$ is flat, it
follows that $y'\in PM\cap N=PN$. But this is false. So $b=st$ for
some $s\in R$. We obtain $t(x+sy'-scm)=0$. Since $M$ is flat and
$tsc\ne 0$ there exists $z\in M$ such that $x=s(-y'+cz)$.

Let $y\in M$. There exists $x\in M\setminus PM$ such that $y\in Rx$. We may assume that $x+PM$ is an element of a basis $(x_i+PM)_{i\in I}$ of $M/PM$. Then $Rx$ is a summand of the free pure submodule $N$ generated by the family $(x_i)_{i\in I}$.

Assume that $R$ is maximal. Let the notations be as above. By \cite[Theorem XI.4.2]{FuSa85} each uniserial $R$-module is pure-injective. So,  $Rx$ is a summand of $M$. Let $u$ be the composition of a projection from $M$ onto $Rx$ with the isomorphism between $Rx$ and $R$. Thus $u\in M^*$ and $u(x)=1$. It follows that $y=u(y)x$. Hence $M$ is locally projective by \cite[Theorem 3.2]{Gar76} or \cite[Theorem 2.1]{ZH76}. \qed

\begin{proposition} \label{P:sinpro}

Let $R$ be a valuation ring such that $Z=P\ne 0$ and let $M$ be a non-zero
$R$-module. Then the following conditions are equivalent:
\begin{enumerate}
\item $M$ is singly projective;
\item $M$ is a flat content module;
\item $M$ is flat and  contains an essential free submodule.
\end{enumerate}
\end{proposition}
\textbf{Proof.} 
$(2)\Rightarrow (1)$ by Proposition~\ref{P:pflat}.

$(1)\Rightarrow (2)$. It remains to show that $M$ is a content module. Let $x\in M$. There exists $y\in M$ and $a\in R$ such that $x=ay$ and $(0:y)=0$. Since  $Z=P$ then $y\notin PM$. We deduce that $\mathrm{c}(x)=Ra$.

$(2)\Leftrightarrow (3)$. By Proposition~\ref{P:content} it remains to show that $(2)\Rightarrow (3)$. Let $N$ be a pure-essential free submodule of $M$. Since $R$ is self FP-injective by \cite[Lemma]{Gil71}, it follows that $N$ is a pure
submodule of each overmodule. So, if $K$ is a submodule of $M$ such
that $K\cap N=0$, then $N$ is a pure submodule of $M/K$, whence $K=0$. \qed 

\begin{corollary} \label{C:selffpinj}
Let $R$ be a valuation ring. The following conditions are equivalent:
\begin{enumerate}
\item $Z=P$;
\item Each singly projective module is a content module.
\end{enumerate} 
\end{corollary}
\textbf{Proof.} It remains to show that $(2)\Rightarrow (1)$. By Proposition~\ref{P:locali} $Q$ is finitely projective over $R$. If $R\ne Q$, then $Q$ is not content on $R$ because, $\forall x\in Q\setminus Z,\ \mathrm{c}(x)=Z$. So $Z=P$. \qed

\begin{corollary} \label{C:hullsinpro} 
Let $R$ be a valuation ring. Then the injective hull of any singly projective module  is 
  singly projective too.
\end{corollary}
\textbf{Proof.} 
 Let $N$ be a non-zero singly projective module. We denote by $E$ its injective hull. For each $s\in
R\setminus Z$ the multiplication by $s$ in $N$ is injective, so the
multiplication by $s$ in $E$ is bijective. Hence $E$ is a
$Q$-module which is flat by Proposition~\ref{P:hullflat}. It is
an essential extension of $N_Z$. From  Propositions~\ref{P:sinpro} and \ref{P:locali}(2) we deduce that $E$ is singly projective.
\qed

\bigskip
Let $R$ be a valuation ring and let $M$ be a non-zero
$R$-module. We say that $M$ is \textbf{separable} if any finite subset  is contained in a summand which is a finite direct sum of uniserial submodules.
      
\begin{corollary} \label{C:sepa}
Let $R$ be a valuation ring. Then any element of a singly
projective module $M$ is contained in a pure uniserial submodule
$L$. Moreover,  if $R$ is maximal, each singly projective module is separable.
\end{corollary}
\textbf{Proof.} By Proposition~\ref{P:content} any element of $M$ is
contained in a pure cyclic free $Q$-submodule $G$ of $M_Z$. We put $L=M\cap G$. As in proof of Proposition~\ref{P:unis}
we show that $L$ is a pure uniserial submodule of $M$. The first assertion is proved.
  
Since $R$ is maximal $L$ is pure-injective by \cite[Theorem XI.4.2]{FuSa85}. So, $L$ is a summand of $M$. Each summand of $M$ is singly projective. It follows that we can complete the proof by induction on the cardinal of the chosen finite subset of $M$. \qed

\begin{corollary} \label{C:sin=loc}
Let $R$ be a valuation ring. Then the following
conditions are equivalent:
\begin{enumerate} 
\item $R$ is self injective;
\item  Each singly projective module is locally projective;
\item  $Z=P$ and each singly projective module is finitely projective.
\end{enumerate}
\end{corollary}
\textbf{Proof.} $(1)\Rightarrow (2)$ by
Proposition~\ref{P:ringinj}.  

 $(2)\Rightarrow (3)$ follows from \cite[Proposition 5.14(4)]{Gar76} and Corollary~\ref{C:selffpinj}.

$(3)\Rightarrow (1)$.
By way of contradiction
suppose that $R$ is not self injective. Let $E$ be the injective hull
of $R$. By Corollary~\ref{C:hullsinpro} $E$ is singly projective. Let
$x\in E\setminus R$ and $M=R+Rx$. Since $E$ is finitely projective, then
there exist a finitely generated free module $F$, a morphism $\phi:M\rightarrow F$ and a
morphism $\pi:F\rightarrow E$ such that $(\pi\circ\phi)(y)=y$ for each
$y\in M$. Let $\tilde{\phi}:M/R\rightarrow F/\phi(R)$ and
$\tilde{\pi}:F/\phi(R)\rightarrow E/R$ be the morphisms induced by
$\phi$ and $\pi$. Then $(\tilde{\pi}\circ\tilde{\phi})(y+R)=y+R$ for
each $y\in M$. Since $\phi(R)$ is a pure submodule of $F$, then
$F/\phi(R)$ is a finitely generated flat module. Hence $F/\phi(R)$ is
free and $E/R$ is singly projective. But $E/R=P(E/R)$ by \cite[Lemma
12]{Couch03}. This contradicts that $E/R$ is a flat content module. By Proposition~\ref{P:sinpro} we conclude that $E=R$. \qed

\begin{corollary} \label{C:sin=fin}
Let $R$ be a valuation ring. Then $Q$ is self injective if and
only if each singly projective module is finitely projective.
\end{corollary}
\textbf{Proof.} By \cite[Theorem 2.3]{KlLe69} $Q$ is maximal if and only if it is self injective. Suppose that $Q$ is self injective and let $M$ be a singly
projective $R$-module. Then $M_Z$ is locally projective over $Q$ by Proposition~\ref{P:locali}(1) and corollary~\ref{C:sin=loc}. Consequently $M$ is finitely projective by Lemma~\ref{L:flat-fproj}. 

Conversely let $M$ be a singly projective $Q$-module. Then $M$ is
singly projective over $R$, whence $M$ is finitely projective over
$R$. If follows that $M$ is finitely projective over $Q$. From
Corollary~\ref{C:sin=loc} we deduce that $Q$ is self injective. \qed

\begin{theorem} \label{T:sepa}
Let $R$ be a valuation ring. Then the following conditions are equivalent:
\begin{enumerate}
\item $R$ is maximal; 
\item each singly projective $R$-module is separable;
\item each flat content module is locally projective.
\end{enumerate} 
\end{theorem}
\textbf{Proof.} $(1)\Rightarrow (2)$ by Corollary~\ref{C:sepa} and $(1)\Rightarrow (3)$ by Proposition~\ref{P:content}.

$(2)\Rightarrow (1):$ let $\widehat{R}$ be
the pure-injective hull of $R$.  By
\cite[Proposition 1 and 2]{Couc06} $\widehat{R}$ is a flat content
module. Consequently $1$  belongs to a  summand $L$ of
$\widehat{R}$ which is a finite direct sum of uniserial modules. But, by \cite[Proposition 5.3]{Fac87} $\widehat{R}$ is
indecomposable. Hence $\widehat{R}$ is uniserial. Suppose that
$R\ne\widehat{R}$. Let $x\in\widehat{R}\setminus R$. Then there exists $c\in P$ such that $1=cx$. Since $R$ is pure in $\widehat{R}$  we get that $1\in P$ which is absurd. Consequently, $R$ is a pure-injective module. So, $R$ is maximal by \cite[Proposition 9]{War69}.

$(3)\Rightarrow (1):$ since $\widehat{R}$ is locally projective then $R$ is a summand of $\widehat{R}$ which is indecomposable. So $R$ is maximal.
\qed

\bigskip
A submodule $N$ of a module $M$ is said to be \textbf{strongly pure} if, $\forall x\in N$ there exists an homomorphism $u:M\rightarrow N$ such that $u(x)=x$. Moreover, if $N$ is pure-essential, we say that $M$ is a \textbf{strongly pure-essential extension} of $N$.

\begin{proposition} Let $R$ be a valuation ring and let $M$ be a flat $R$-module. Then $M$
  is locally
  projective if and only if it is a strongly pure-essential extension of a free
  module.
\end{proposition}
\textbf{Proof.} Let $M$ be a non-zero locally projective
$R$-module. Then $M$ is a flat content module. So $M$ contains a
pure-essential free submodule $N$.  Let $x\in N$. There exist $u_1,\dots ,u_n\in M^*$ and $y_1,\dots ,y_n\in M$ such that $x=\sum_{i=1}^nu_i(x)y_i$. Since $N$ is a pure submodule, $y_1,\dots ,y_n$ can be chosen in $N$. Let $\phi :M\rightarrow N$ be the homomorphism defined by $\phi(z)=\sum_{i=1}^nu_i(z)y_i$. Then $\phi(x)=x$. So, $N$ is a strongly pure submodule of $M$.

Conversely, assume that $M$ is a strongly pure-essential extension of
a free submodule $N$. Let $x\in M$.  As in proof of Proposition~\ref{P:content}, $x=s(y+cz)$,
where $s\in R,\ y\in N\setminus PN,\ c\in P$ and $z\in M$. Since $N$
is strongly pure, there exists a morphism $\phi:M\rightarrow N$ such
that $\phi(y)=y$. Let $\{e_i\mid i\in I\}$ be a basis of $N$. Then $y=\sum_{i\in J}a_ie_i$ where $J$ is a finite subset of $I$ and $a_i\in R$, $\forall i\in J$. Since $y\in N\setminus PN$ there exists $j\in J$ such that $a_j\notin P$. We easily check that $\{y,e_i\mid i\in I,\ i\ne j\}$ is a basis of $N$ too. Hence $Ry$ is a summand of $N$. Let $u$ be the composition of $\phi$ with a projection
of $N$ onto $Ry$ and with the isomorphism between $Ry$ and $R$. Then
$u\in M^*,\ u(y)=1$ and $u(y+cz)=1+cu(z)=v$ is a unit. We put
$m=v^{-1}(y+cz)$. It follows that $x=u(x)m$. Hence $M$ is locally
projective by \cite[Theorem 3.2]{Gar76} or \cite[Theorem 2.1]{ZH76}. \qed

\begin{corollary} Let $R$ be a valuation ring and let $M$ be a locally
  projective $R$-module. If $M/PM$ is finitely generated then $M$ is
  free.
\end{corollary}

\begin{theorem} Let $R$ be a valuation ring. The following conditions are equivalent:
\begin{enumerate}
\item $Z$ is nilpotent;
\item $Q$ is an artinian ring;
\item Each flat $R$-module is finitely projective;
\item Each flat $R$-module is singly projective.
\end{enumerate} 
\end{theorem}
\textbf{Proof.} $(1)\Leftrightarrow (2)$. If $Z$ is nilpotent then $Z^2\ne Z$. It follows that $Z$ is finitely generated over $Q$ and it is the only prime ideal of $Q$. So, $Q$ is artinian. The converse is well known.

$(2)\Rightarrow (3)$ is a consequence of \cite[Corollary 7]{She91} and
it is obvious that $(3)\Rightarrow (4)$.

$(4)\Rightarrow (2)$. First we prove that each flat $Q$-module is singly projective. By Proposition~\ref{P:sinpro} it follows that each flat $Q$-module is content. We deduce that $Q$ is perfect by Theorem~\ref{T:perf}. We conclude that $Q$ is artinian since $Q$ is a valuation ring.
\qed

\section{Strongly coherence or $\pi$-coherence of valuation rings.}
\label{S:coherence}

In this section we study the valuation rings, with non-zero zero-divisors, for which any product of content (respectively singly, finitely, locally projective) modules is content (respectively singly, finitely, locally projective) too.
\begin{theorem} \label{T:pi-cont} 
Let $R$ be a valuation ring such that $Z\ne 0$. Then the following conditions are
  equivalent:
\begin{enumerate}
\item Each product of content modules is content;
\item $R^R$ is a content module;
\item For each ideal $A$ there exists $a\in R$ such that either $A=Ra$
  or $A=Pa$;
\item Each ideal is countably generated and $R^{\mathbb{N}}$ is a content module;
\item The intersection of any non-empty family of principal ideals is finitely generated.
\end{enumerate}
\end{theorem}
\textbf{Proof.} The conditions $(1),\ (2)$ and $(5)$ are equivalent by \cite[Theorem 5.15]{Gar76}. By \cite[Lemma 29]{Couch03} $(3)\Leftrightarrow (5)$.

$(2)\Rightarrow (4)$. It is obvious that $R^{\mathbb{N}}$ is a content module. Since $(2)\Leftrightarrow (3)$ then $P$ is the only prime ideal. We conclude by \cite[Corollary 36]{Couch03}.

$(4)\Rightarrow (3)$. Let $A$ be a non-finitely generated ideal. Let $\{a_n\mid n\in\mathbb{N}\}$ be a spanning set of $A$. Then $x=(a_n)_{n\in\mathbb{N}}\in R^{\mathbb{N}}$. It follows that $x=ay$ for some $a\in \mathrm{c}(x)$ and $y\in R^{\mathbb{N}}$, and $\mathrm{c}(x)=Ra$ . So, if $y=(b_n)_{n\in\mathbb{N}}$, we easily check that $P$ is generated by $\{b_n\mid n\in\mathbb{N}\}$. Hence $A=aP$. \qed

\bigskip

By Proposition~\ref{P:pflat} each valuation domain $R$ verifies the  first two
conditions of the next theorem.

\begin{theorem} \label{T:pi-sp} 
Let $R$ be a valuation ring such that $Z\ne 0$. Then the following conditions are
  equivalent:
\begin{enumerate}
\item Each product of singly projective modules is singly
  projective.
\item $R^R$ is singly projective;
\item $C^*$ is a finitely generated  module for
  each cyclic   module $C$;
\item $(0:A)$ is finitely generated for each proper ideal $A$;
\item $P$ is principal or faithful and for each ideal $A$ there exists $a\in R$ such that either $A=Ra$
  or $A=Pa$;
\item Each ideal is countably generated and $R^{\mathbb{N}}$ is singly projective;
\item Each product of flat content modules is flat content;
\item $R^R$ is a flat content module;
\item Each ideal is countably generated and $R^{\mathbb{N}}$ is a flat content module;
\item $P$ is principal or faithful and the intersection of any non-empty family of principal ideals is finitely generated.
\end{enumerate}
\end{theorem} 
\textbf{Proof.} It is obvious that $(1)\Rightarrow (2)$ and $(7)\Rightarrow (8)$.

$(3)\Leftrightarrow (4)$ because $(0:A)\cong (R/A)^*$.

$(2)\Rightarrow (4)$. Let $A$ be a proper ideal. Then $R^A$ is singly projective too and $x=(a)_{a\in A}$ is an element of
$R^A$. Therefore $x$ belongs to a 
cyclic free submodule of $R^A$ by Lemma~\ref{L:sinval}. Since $R^R$ is flat, $R$ is coherent by \cite[Theorem IV.2.8]{FuSa01}. Consequently $(0:A)=(0:x)$
is finitely generated.

$(4)\Rightarrow (5)$.  Then $R$ is coherent because $R$ is a valuation ring. Since $Z\ne 0$, $Z=P$ by \cite[Theorem 10]{Couch03}. If $P$ is not finitely generated then $P$ cannot
be an annihilator. So $P$ is faithful. By \cite[Lemma 3]{Gil71} and \cite[Proposition 1.3]{KlLe69}, if $A$ is a proper ideal then
either $A=(0:(0:A))$ or $A=P(0:(0:A))$. By $(4)$, $(0:(0:A))=Ra$ for
some $a\in P$.

$(5)\Rightarrow (1)$. Let $(M_i)_{i\in I}$ be a family of singly
projective modules. Let $x=(x_i)_{i\in I}$ be an element of $\Pi_{i\in
  I}M_i$. Since $M_i$ is singly projective for each $i\in I$ there
exist $a_i\in R$ and $y_i\in M_i$ such that $x_i=a_iy_i$ and
$(0:y_i)=0$. We have either $\sum_{i\in I}Ra_i=Ra$ or $\sum_{i\in
  I}Ra_i=Pa$ for some $a\in R$. Then, $\forall i\in I,\ \exists b_i\in
R$ such that $a_i=ab_i$. Therefore either $\exists i\in I$ such that
$b_i$ is a unit, or $P=\sum_{i\in I}Rb_i$. It follows that $x=ay$ where
$y=(b_iy_i)_{i\in I}$. Now it is easy to check that $(0:y)=0$.

$(6)\Rightarrow (4)$. Since each ideal is countably generated then so is each submodule of a finitely generated free module. So, the flatness of $R^\mathbb{N}$ implies that $R$ is coherent. Let $A$ be a proper ideal generated by $\{a_n\mid n\in\mathbb{N}\}$. Then $x=(a_n)_{n\in\mathbb{N}}$ is an element of $R^\mathbb{N}$. Therefore $x$ belongs to a cyclic free submodule of $R^\mathbb{N}$ by Lemma~\ref{L:sinval}. Consequently $(0:A)=(0:x)$
is finitely generated because $R$ is coherent.

$(5)\Rightarrow (9)$. By Theorem~\ref{T:pi-cont}($(3)\Leftrightarrow (4)$) it remains to show that $R^\mathbb{N}$ is flat. This is true because  $(5)\Rightarrow (1)$. 

$(1)\Leftrightarrow (7)$. Since $(1)\Rightarrow (2)$ or $(7)\Rightarrow (8)$, $R$ is coherent. From $Z\ne 0$ and \cite[Theorem 10]{Couch03} it follows that $Z=P$. Now we use Proposition~\ref{P:sinpro} to conclude.

$(2)\Leftrightarrow (8)$. Since $R^R$ is flat, $R$ is coherent. We do as above to conclude.

$(6)\Leftrightarrow (9)$. Since each submodule of a free module of finite rank is countably generated, then the flatness of $R^{\mathbb{N}}$ implies that $R$ is coherent. So we conclude as above.

$(5)\Leftrightarrow (10)$ by Theorem~\ref{T:pi-cont}($(3)\Leftrightarrow (5)$). 

The last assertion is already shown. So, the proof is complete.
\qed

\begin{remark}
\textnormal{When $R$ is a valuation domain, the conditions $(5),\ (7),\ (8),\ (9)$ and $(10)$ are equivalent by \cite[Theorem 4]{Kem00} and \cite[Corollary 36]{Couch03}.}
\end{remark}

\begin{remark}
\textnormal{If $R$ satisfies  the conditions of Theorem~\ref{T:pi-cont} and if $P$ is not faithful and not finitely generated then $R$ is not coherent and doesn't satisfy the conditions of Theorem~\ref{T:pi-sp}.}
\end{remark}

\bigskip 
By \cite[Corollary 3.5]{Fra84} or \cite[Theorem 3]{Kem00}, a valuation domain $R$ is
strongly coherent if and only if either its order group is $\mathbb{Z}$ or if
$R$ is maximal and its order group is $\mathbb{R}$. It is easy to
check that each Pr\"ufer domain is $\pi$-coherent because it satisfies the fourth condition of the
next theorem. When $R$ is a valuation ring with non-zero zero-divisors
we get:

\begin{theorem} \label{T:pi-coh} 
Let $R$ be a valuation ring such that $Z\ne 0$. Then the following conditions are
  equivalent:
\begin{enumerate}
\item $R$ is strongly coherent;
\item $R$ is $\pi$-coherent;
\item $R^R$ is singly projective and separable;
\item $C^*$ is a finitely generated  module for
  each finitely generated   module $C$;
\item $(0:A)$ is finitely generated for each proper ideal
  $A$ and $R$ is self injective;
\item $R$ is maximal, $P$ is principal or faithful and for each ideal $A$ there exists $a\in R$ such that either $A=Ra$
  or $A=Pa$;
\item Each ideal is countably generated and $R^{\mathbb{N}}$ is singly projective and separable;
\item $R^R$ is a separable flat content module;
\item Each ideal is countably generated and $R^{\mathbb{N}}$ is a separable flat content module;
\item Each product of separable flat content modules is a separable flat content module;
\item $R$ is maximal, $P$ is principal or faithful and the intersection of any non-empty family of principal ideals is finitely generated.
\end{enumerate}
\end{theorem} 
\textbf{Proof.}
 By Theorem~\ref{T:picoh} $(1)\Rightarrow (2)$. It is obvious that $(10)\Rightarrow (8)$.  By \cite[Theorem 1]{Cam90}
$(2)\Leftrightarrow (4)$. By Theorem~\ref{T:pi-sp}, Theorem~\ref{T:sepa} and \cite[Theorem 2.3]{KlLe69}
 $(5)\Leftrightarrow (6)$ and $(6)\Rightarrow (7)$ . By Theorem~\ref{T:pi-sp} $(6)\Leftrightarrow (11)$, $(7)\Leftrightarrow (9)$ and $(3)\Leftrightarrow (8)$.
 
$(4)\Rightarrow (6)$. By Theorem~\ref{T:pi-sp} $R$ is coherent and self
FP-injective and it remains to prove that $R$ is
maximal if $P$ is not principal.
Let $E$ be the injective hull of $R$. If $R\ne E$ let $x\in E\setminus
R$. Since $R$ is an essential submodule of $E$, $(R:x)=rP$ for some
$r\in R$. Then $(R:rx)=P$. Let $M$ be the submodule of $E$ generated
by $1$ and $rx$. We put $N=M/R$. Then $N\cong R/P$. We get that
$N^*=0$ and $M^*$ is isomorphic to a principal
ideal of $R$. Moreover, since $(R:rx)=P$, for each $t\in P$ the
multiplication by $t$ in $M$ is a non-zero element of
$M^*$. Since $P$ is faithful we get that
$M^*\cong R$. Let $g\in M^*$ such that
the restriction of $g$ to $R$ is the identity. For each $p\in P$ we
have $pg(rx)=prx$. So $(0:g(rx)-rx)=P$. Since $P$ is faithful, there
is no simple submodule in $E$. Hence $g(rx)=rx$ but this is not possible
because $g(rx)\in R$ and $rx\notin R$. Consequently $R$ is
self-injective and  maximal.

$(2)\Rightarrow (1)$. Since $(2)\Rightarrow (6)$ $R$ is self injective. We conclude by proposition~\ref{P:ringinj}.

$(3)\Rightarrow (1)$. Since $R^R$ is singly projective, by Theorem~\ref{T:pi-sp} $R$ is coherent and self FP-injective. So, if $U$ is a uniserial summand of $R^R$, then $U$ is singly projective and consequently $U\ne PU$. Let $x\in U\setminus PU$. It is easy to check that $U=Rx$ and that $(0:x)=0$. Hence $R^R$ is locally projective and $R$ is strongly coherent.

$(7)\Rightarrow (4)$.  Let $F_1\rightarrow F_0\rightarrow C\rightarrow 0$ be a free presentation of a finitely generated module $C$, where $F_0$ is finitely generated. It follows that $F_1$ is countably generated. As above we prove that $R^{\mathbb{N}}$ is locally projective. By Theorem~\ref{T:pi-sp} $R$ is coherent and consequently each finitely generated submodule of $R^{\mathbb{N}}$ is finitely presented. Since $F_1^*\cong R^{\mathbb{N}}$ we easily deduce that $C^*$ is finitely generated. 

$(1)\Rightarrow (10)$. Since $(1)\Rightarrow (6)$, $R$ is maximal. We use Theorem~\ref{T:sepa} to conclude. The proof is now complete.
\qed


\begin{thebibliography}{10}

\bibitem{Azu87}
G.~Azumaya.
\newblock Finite splitness and finite projectivity.
\newblock {\em J. Algebra}, 106:114--134, (1987).

\bibitem{Bas60}
H.~Bass.
\newblock Finistic dimension and a homological generalization of semi-primary
  rings.
\newblock {\em Trans. Amer. Math. Soc.}, 95:466--488, (1960).

\bibitem{Cam90}
V.~Camillo.
\newblock Coherence for polynomial rings.
\newblock {\em J. Algebra}, 132:72--76, (1990).

\bibitem{Couch03}
F.~Couchot.
\newblock Injective modules and fp-injective modules over valuations rings.
\newblock {\em J. Algebra}, 267:359--376, (2003).

\bibitem{Couc06}
F.~Couchot.
\newblock Pure-injective hulls of modules over valuation rings.
\newblock {\em J. Pure Appl. Algebra}, 207:63--76, (2006).

\bibitem{CoPe70}
H.~Cox and R.~Pendleton.
\newblock Rings for which certain flat modules are projective.
\newblock {\em Trans. Amer. Math. Soc.}, 150:139--156, (1970).

\bibitem{Fac87}
A.~Facchini.
\newblock Relative injectivity and pure-injective modules over {P}r\"ufer
  rings.
\newblock {\em J. Algebra}, 110:380--406, (1987).

\bibitem{Fra84}
B.~Franzen.
\newblock On the separability of a direct product of free modules over a
  valuation domain.
\newblock {\em Arch. Math.}, 42:131--135, (1984).

\bibitem{FuSa85}
L.~Fuchs and L.~Salce.
\newblock {\em Modules over valuation domains}, volume~97 of {\em Lecture Notes
  in Pure an Appl. Math.}
\newblock Marcel Dekker, New York, (1985).

\bibitem{FuSa01}
L.~Fuchs and L.~Salce.
\newblock {\em Modules over Non-Noetherian Domains}.
\newblock Number~84 in Mathematical Surveys and Monographs. American
  Mathematical Society, Providence, (2001).

\bibitem{Gar76}
G.S. Garfinkel.
\newblock Universally {T}orsionless and {T}race {M}odules.
\newblock {\em Trans. Amer. Math. Soc.}, 215:119--144, (1976).

\bibitem{Gil71}
D.T. Gill.
\newblock Almost maximal valuation rings.
\newblock {\em J. London Math. Soc.}, 4:140--146, (1971).

\bibitem{GrRa71}
L.~Gruson and M.~Raynaud.
\newblock Crit\`ers de platitude et de projectivit\'e.
\newblock {\em Invent. Math.}, 13:1--89, (1971).

\bibitem{Ja73}
S.~Jain.
\newblock Flatness and {FP}-injectivity.
\newblock {\em Proc. Amer. Math. Soc.}, 41(2):437--442, (1973).

\bibitem{Jon82}
M.F. Jones.
\newblock Flatness and f-projectivity of torsion-free modules and injective
  modules.
\newblock In {\em Advances in Non-commutative Ring Theory}, number 951 in
  Lecture Notes in Mathematics, pages 94--116, New York/Berlin, (1982).
  Springer-Verlag.

\bibitem{Kem00}
R.~Kemper.
\newblock Product-trace-rings and a question of {G.} {S.} {G}arfinkel.
\newblock {\em Proc. Amer. Math. Soc.}, 128(3):709--712, (2000).

\bibitem{KlLe69}
G.B. Klatt and L.S. Levy.
\newblock Pre-self injectives rings.
\newblock {\em Trans. Amer. Math. Soc.}, 137:407--419, (1969).

\bibitem{Kob84}
S.~Kobayashi.
\newblock A note on regular self-injective rings.
\newblock {\em Osaka J. Math.}, 21(3):679--682, (1984).

\bibitem{RaRa73}
V.S. Ramamurthi and K.M. Rangaswamy.
\newblock On finitely injective modules.
\newblock {\em J. Aust. Math. Soc.}, XVI(2):239--248, (1973).

\bibitem{She91}
Zhu Shenglin.
\newblock On rings over which every flat left module is finitely projective.
\newblock {\em J. Algebra}, 139:311--321, (1991).

\bibitem{War69}
R.B. Warfield.
\newblock Purity and algebraic compactness for modules.
\newblock {\em Pac. J. Math.}, 28(3):689--719, (1969).

\bibitem{ZH76}
B.~Zimmermann-Huisgen.
\newblock Pure submodules of direct products of free modules.
\newblock {\em Math. Ann.}, 224:233--245, (1976).

\bibitem{ZH80}
B.~Zimmermann-Huisgen.
\newblock Direct products of modules and algebraic compactness.
\newblock Habilitationsschrift, Tech. Univ. M\"unchen, (1980).

\end{thebibliography}
\end{document}